# STOCHASTIC BOUNDS FOR LÉVY PROCESSES

By R. A. Doney

*University of Manchester*

Using the Wiener–Hopf factorization, it is shown that it is possible to bound the path of an arbitrary Lévy process above and below by the paths of two random walks. These walks have the same step distribution, but different random starting points. In principle, this allows one to deduce Lévy process versions of many known results about the large-time behavior of random walks. This is illustrated by establishing a comprehensive theorem about Lévy processes which converge to $\infty$ in probability.

**1. Introduction.** If $X = (X_t, t \geq 0)$ is an arbitrary Lévy process, we would frequently like to be able to assert that some aspect of its behavior as $t \to \infty$ can be seen to be true "by analogy with known results for random walks." An obvious way to try to justify such a claim is via the random walk $S^{(\delta)} := (X(n\delta), n \geq 0)$, for fixed $\delta > 0$. (This process is often called the $\delta$-skeleton of $X$.) However, it is often difficult to control the deviation of $X$ from $S^{(\delta)}$. A further problem stems from the fact that the distribution of $S_1^{(\delta)} = X(\delta)$ is determined via the Lévy–Khintchine formula and not directly in terms of the characteristics of $X$, that is, the Lévy measure $\Pi$, the Brownian coefficient $\sigma^2$, and $\gamma$, the linear coefficient in the Lévy–Khintchine formula.

An alternative approach is to use the random walk which results from observing $X$ at the times at which its "large jumps" occur. Specifically, we will assume, here and throughout, that $\Pi(\mathbb{R}) > 0$, since otherwise $X$ is a Brownian motion or a pure drift and all the results we give are already known. Then we take a fixed interval $I = [-\eta_1, \eta_2]$ which contains zero and has $\Delta := \Pi(I^c) > 0$, put $\tau_0 = 0$, and for $n \geq 1$ write $\tau_n$ for the time at which $J_n$, the $n$th jump in $X$ whose value lies in $I^c$, occurs. (It might seem to be natural to assume that $\eta_1 = \eta_2$, but it is possible that the extra generality can be useful.) The random walk is then defined by

$$\hat{S} := (\hat{S}_n, n \geq 0), \qquad \text{where } \hat{S}_n = X(\tau_n). \tag{1.1}$$

Of course $(\tau_n, n \geq 1)$ is a Poisson process of rate $\Delta$ which is independent of $(J_n, n \geq 1)$, and this latter is a sequence of i.i.d. (independent, identically distributed) random variables having the distribution $\Delta^{-1}\mathbf{1}_{I^c}\Pi(dx)$. We will write $\hat{Y}_1, \hat{Y}_2, \ldots$ for the steps in $\hat{S}$, so that, with $e_r := \tau_r - \tau_{r-1}$ and $r \geq 1$,

$$\hat{Y}_r = X(\tau_r) - X(\tau_{r-1}) = J_r + \tilde{X}(\tau_r) - \tilde{X}(\tau_{r-1}) \stackrel{D}{=} J_r + \tilde{X}(e_r), \tag{1.2}$$

where $\tilde{X}$ is "$X$ with the jumps $J_1, J_2, \ldots$ removed." This is also a Lévy process whose Lévy measure is the restriction of $\Pi$ to $I$. Furthermore, $\tilde{X}$ is independent of $\{(J_n, \tau_n), n \geq 1\}$, and since it has no large jumps, it follows that $E\{e^{\lambda \tilde{X}_t}\}$ is finite for all real $\lambda$. Thus the contribution of $\sum_1^n \tilde{X}(e_r)$ to $\hat{S}_n$ can be easily estimated, and for many purposes $\hat{Y}_r$ can be replaced by $J_r + \tilde{\mu}$, where $\tilde{\mu} = E\tilde{X}(\tau_1)$. In order to control the deviation of $X$ from $\hat{S}$, it is natural to use the stochastic bounds

$$I_n \leq X_t \leq M_n \qquad \text{for } \tau_n \leq t < \tau_{n+1}, \tag{1.3}$$

where

$$I_n := \inf_{\tau_n \leq t < \tau_{n+1}} X_t, \qquad M_n := \sup_{\tau_n \leq t < \tau_{n+1}} X_t, \tag{1.4}$$

and write

$$M_n = \hat{S}_n + \tilde{m}_n \quad \text{and} \quad I_n = \hat{S}_n + \tilde{i}_n. \tag{1.5}$$

Here

$$\tilde{m}_n = \sup_{0 \leq s < e_{n+1}} \{\tilde{X}(\tau_n + s) - \tilde{X}(\tau_n)\}, \qquad n \geqslant 1, \tag{1.6}$$

$$\tilde{i}_n = \inf_{0 \leq s < e_{n+1}} \{\tilde{X}(\tau_n + s) - \tilde{X}(\tau_n)\}, \qquad n \geqslant 1, \tag{1.7}$$

are each i.i.d. sequences, and both $\tilde{m}_n$ and $\tilde{i}_n$ are independent of $\hat{S}_n$. This method also leads to some technical complications; see, for example, the proofs of Theorems 3.3 and 3.4 in [2].

The main point of this paper is to demonstrate that there is a different representation for the random variables $M_n$ and $I_n$ in (1.4) which allows us to draw conclusions about the asymptotic behavior of Lévy processes from the corresponding results for random walks in a simpler way.

THEOREM 1.1. *Using the above notation we have, for any fixed $\eta_1, \eta_2 > 0$ with $\Delta = \Pi(I^c) > 0$,*

$$(1.8) \qquad M_n = S_n^{(+)} + \tilde{m}_0, \qquad I_n = S_n^{(-)} + \tilde{i}_0, \qquad n \geq 0,$$

*where both of the processes $S^{(+)} = (S_n^{(+)}, n \geq 0)$ and $S^{(-)} = (S_n^{(-)}, n \geq 0)$ are random walks with the same distribution as $\hat{S}$. Moreover, $S^{(+)}$ and $\tilde{m}_0$ are independent, as are $S^{(-)}$ and $\tilde{i}_0$.*

Comparing the representations (1.5) and (1.8), note that, for each fixed $n$, the pairs $(\hat{S}_n, \tilde{m}_n)$ and $(S_n^{(+)}, \tilde{m}_0)$ have the same joint law; however, the latter representation has the great advantage that the term $\tilde{m}_0$ does not depend on $n$.

A straightforward consequence of Theorem 1.1 is

PROPOSITION 1.2. *Suppose that $b \in RV(\alpha)$, and $\alpha > 0$. Then, for any fixed $\eta_1, \eta_2 > 0$ with $\Delta = \Pi(I^c) > 0$, and any $c \in [-\infty, \infty]$,*

$$(1.9) \quad \frac{\hat{S}_n}{b(n)} \xrightarrow{a.s.} c \quad \text{as } n \to \infty \quad \Longleftrightarrow \quad \frac{X_t}{b(t)} \xrightarrow{a.s.} \frac{c}{\Delta^\alpha} \quad \text{as } t \to \infty.$$

[Here $RV(\alpha)$ denotes the class of positive functions which are regularly varying with index $\alpha$ at $\infty$.]

From this, and analogous statements for lim sup and lim inf, known results about Lévy processes such as strong laws and laws of the iterated logarithm can easily be deduced. But there is a vast literature on the asymptotic behavior of random walks, and by no means all the results it contains have been extended to the setting of Lévy processes. Using Theorem 1.1 we can show, for example, that the classical results of [6] about strong limit points of random walks, and results about the lim sup behavior of $S_n/n^\alpha$ and $|S_n|/n^\alpha$ and hence about first passage times outside power-law type boundaries in [11], all carry over easily. A further case in point is that of existence of moments for first and last passage times in the transient case; see [5] and [8], which completed results of many earlier authors. It turns out that the combination of the stochastic bound (1.3) and Theorem 1.1 is ideally suited to analyzing the corresponding Lévy situation; see [3].

However, here we will concentrate on the extensive results which have emerged in a series of papers by Kesten and Maller which explore various aspects of the asymptotic behavior of random walks which converge to $+\infty$ *in probability*. The following theorem gives Lévy process versions of just a sample of their results; specifically Theorem 2.1 in [7], Theorem 3 in [9] and Proposition 1 in [10].

We will use the following notation: our Lévy process will be written as

$$(1.10) \qquad X_t = \gamma t + \sigma B_t + Y_t^{(1)} + Y_t^{(2)},$$

where $B$ is a standard BM, $Y^{(1)}$ is a pure jump martingale formed from the jumps whose absolute values are less than or equal to 1, $Y^{(2)}$ is a compound Poisson process formed from the jumps whose absolute values exceed 1, and $B$, $Y^{(1)}$ and $Y^{(2)}$ are independent. We will assume throughout that $\Pi(\mathbb{R}) > 0$.

For $x > 0$, we introduce the tail functions

$$(1.11) \qquad N(x) = \Pi\{(x, \infty)\}, \qquad M(x) = \Pi\{(-\infty, -x)\},$$

and the tail sum and difference

$$(1.12) \quad T(x) = N(x) + M(x), \qquad D(x) = N(x) - M(x), \qquad x > 0.$$

The rôles of truncated first and second moments are played by

$$(1.13) \qquad \begin{aligned} A(x) &= \gamma + D(1) + \int_1^x D(y)\,dy, \\ U(x) &= \sigma^2 + 2\int_0^x y T(y)\,dy, \qquad x > 0. \end{aligned}$$

Both $A$ and $U$ are continuous functions with $A(x)/x \to 0$ and $U(x)/x^2 \to 0$ as $x \to \infty$. Finally, we introduce the two-sided exit time by

$$(1.14) \qquad T_r = \inf\{t : |X_t| > r\}.$$

THEOREM 1.3. *Assume $M(x) > 0$ for all $x > 0$. Then the following are equivalent:*

$$(1.15) \qquad P(X_{T_r} > 0) \to 1 \qquad \text{as } r \to \infty;$$

$$(1.16) \qquad P(X_t > 0) \to 1 \qquad \text{as } t \to \infty;$$

$$(1.17) \qquad X_t \xrightarrow{P} +\infty \qquad \text{as } t \to \infty;$$

$$(1.18) \qquad \frac{X_t}{b(t)} \xrightarrow{P} +\infty \qquad \text{as } t \to \infty \text{ for some } b \in RV(1);$$

$$(1.19) \qquad \frac{A(x)}{\sqrt{U(x)M(x)}} \to +\infty \qquad \text{as } x \to \infty.$$

REMARK 1.1. The assumption that $M(x) > 0$ for all $x > 0$ is not essential; in the contrary case, the theorem still holds, except that (1.19) should be replaced by the condition $EX_1 > 0$. (Note that since $EX_1^- < \infty$ in this case, $EX_1$ is well defined.)

REMARK 1.2. On the basis of Theorem 1 and Remark (viii) of [9], it is natural to suppose that the following subsequential version of Theorem 1.3 is valid:

THEOREM 1.4. *Assume $M(x) > 0$ for all $x > 0$. Then the following are equivalent:* (1.15) *holds for some (deterministic) sequence $r_k \to \infty$;* (1.16) *holds for some (deterministic) sequence $t_k \to \infty$;* (1.17) *holds for some (deterministic) sequence $t_k \to \infty$;* (1.19) *holds for some (deterministic) sequence $x_k \to \infty$; and, for some (deterministic) sequence $t_k \to \infty$,*

$$\frac{X_{t_k}}{\sqrt{t_k}} \xrightarrow{P} +\infty \qquad as\ t \to \infty.$$

This is in fact correct, and can be proved by arguments that are similar to, but more complicated than, those we use to prove Theorem 1.3, but we omit the details.

REMARK 1.3. It is not difficult to see that the results established in [10] about random walks which leave regions of the form $\{(x, n) : |x| \leq rn^\kappa\}$ at the upper boundary with probability approaching 1 as $r \to \infty$ can also be shifted to the Lévy process setting by analogous arguments.

REMARK 1.4. In the applications discussed here, we work with a fixed choice of the cut-off points $\eta_1$ and $\eta_2$. However, provided $\Pi(\mathbb{R}) = \infty$, we could get a sequence of bounds by taking $\eta_i^{(n)} \downarrow 0$ as $n \to \infty$, $i = 1, 2$. It is not difficult to see that such a sequence would converge uniformly to $X$ a.s. on compact time intervals. This fact might have other applications, for example in the important area of simulation.

## 2. Proofs.

PROOF OF THEOREM 1.1. The Wiener–Hopf factorization for $\tilde{X}$ (see [1], page 165) asserts that the random variables $\tilde{m}_0 = \sup_{0 \leq t < e_1} \tilde{X}_t$ and $\tilde{X}_{e_1} - \tilde{m}_0$ are independent, and that the latter has the same distribution as $\tilde{i}_0 = \inf_{0 \leq t < e_1} \tilde{X}_t$. [Recall that $\tilde{X}$ and $e_1$ are independent and $e_1$ has an $\text{Exp}(\Delta)$ distribution.] Since

$$\begin{aligned}
M_1 &= \sup_{e_1 \leq t < e_1 + e_2} X_t = \tilde{X}(e_1) + J_1 + \sup_{0 \leq t < e_2}\{\tilde{X}(e_1 + t) - \tilde{X}(e_1)\} \\
&= \tilde{m}_0 + \{\tilde{X}(e_1) - \tilde{m}_0\} + J_1 + \tilde{m}_1 \\
&:= \tilde{m}_0 + Y_1^{(+)},
\end{aligned}$$

where all four random variables in the second line are independent, we see that $Y_1^{(+)}$ is independent of $\tilde{m}_0$ and has the same distribution as $J_1 + \tilde{X}(e_1)$, and hence as $X(e_1)$. A similar calculation applied to $M_n$ gives the required conclusions for $S^{(+)}$, and since $S^{(-)}$ is $S^{(+)}$ evaluated for $-X$, the proof is completed. $\square$

PROOF OF PROPOSITION 1.2. With $N_t = \max\{n : \tau_n \leq t\}$, we have, from (1.3) and (1.8),

$$(2.1) \qquad \frac{\tilde{i}_0}{b(t)} + \frac{S^{(-)}_{N_t}}{b(N_t)} \cdot \frac{b(N_t)}{b(t)} \leq \frac{X_t}{b(t)} \leq \frac{S^{(+)}_{N_t}}{b(N_t)} \cdot \frac{b(N_t)}{b(t)} + \frac{\tilde{m}_0}{b(t)}.$$

Clearly the extreme terms converge a.s. to zero as $t \to \infty$, and by the strong law $b(N_t)/b(t) \stackrel{\text{a.s.}}{\to} 1/\Delta^\alpha$. So if $\frac{\hat{S}_n}{b(n)} \stackrel{\text{a.s.}}{\to} c$ as $n \to \infty$, then $\frac{S^{(+)}_n}{b(n)} \stackrel{\text{a.s.}}{\to} c$ and $\frac{S^{(-)}_n}{b(n)} \stackrel{\text{a.s.}}{\to} c$, and hence $\frac{X_t}{b(t)} \stackrel{\text{a.s.}}{\to} \frac{c}{\Delta^\alpha}$ as $t \to \infty$. On the other hand, if this last is true, we can use (2.1) with $t = \tau_n$ to reverse the argument. $\square$

PROOF OF THEOREM 1.3. Since $M(1) > 0$ by assumption, throughout this proof we will take $\eta_1 = \eta_2 = 1$; note that this means that, in the notation of (1.10),

$$(2.2) \quad \tilde{X}_t = \gamma t + B_t + Y_t^{(1)} \quad \text{and} \quad Y_t^{(2)} = \hat{S}_{N_t}, \qquad \text{where } \hat{S}_n = \sum_1^n \hat{Y}_m.$$

Consequently, we have $E\tilde{X}_1 = \gamma$. We will first use Theorem 1.1 to show that (1.15) is equivalent to (1.19). Recall that Theorem 3 of [9] states, for any random walk $S$, that $P(S_n > 0) \to 1$ if and only if $P(S$ exits $(-r, r)$ at the top$) \to 1$. Next note that it is easy to see that this statement still holds if we replace the interval $(-r, r)$ by $(-r + b, r + c)$ for any fixed $b$ and $c$. It also follows easily from (1.3) and (1.8) that, for any fixed $a > 0$,

$\limsup P(X_{T_r} > 0)$

$\qquad \leq \limsup P\{(M_n, n \geq 0) \text{ exits } (-r - a, r + a) \text{ at the top}\}$

$\qquad \leq \limsup \int P\{S^{(+)} \text{ exits } (-r - a + x, r + a + x) \text{ at the top}\} P(\tilde{m}_0 \in dx).$

From this and a corresponding inequality for $S^{(-)}$, we conclude that

$$(2.3) \qquad\qquad (1.15) \quad \Longleftrightarrow \quad P(\hat{S}_n > 0) \to 1.$$

Now note that if we write $\hat{Y}_m = J_m + \gamma/\Delta + Z_m$, so that $Z_m = \tilde{X}(e_m) - \gamma/\Delta$, we have

$$EZ_1 = E(e_1)E\tilde{X}_1 - \gamma/\Delta = 0.$$

Since $\operatorname{Var} Z_1 < \infty$, it is clear that if $b \in RV(1)$, then $\sum_1^n Z_m/b_n \stackrel{P}{\to} 0$; thus $\hat{S}_n/b_n \stackrel{P}{\to} \infty$ if and only if $S^*_n/b_n \stackrel{P}{\to} \infty$, where $S^*_n = \sum_1^n J^*_m$ and $J^*_m = J_m +$ $\gamma/\Delta$. Since the condition $P(S_n > 0) \to 1$ is equivalent, for any random walk, to the existence of such a $b$ with $S_n/b_n \xrightarrow{P} \infty$ (see Proposition 1 of [10]), another appeal to Theorem 3 of [9] shows we have established that

(2.4)
$$(1.15) \iff P(S_n^* > 0) \to 1$$
$$\iff \frac{A^*(x)}{\sqrt{U^*(x)F^*(-x)}} \to +\infty \quad \text{as } x \to \infty,$$

where $F^*$ is the distribution function of $J_1^*$ and

$$A^*(x) = \int_0^x \{1 - F^*(y) - F^*(-y)\}\, dy,$$

$$U^*(x) = 2\int_0^x y\{1 - F^*(y) - F^*(-y)\}\, dy.$$

However, since $F^*(dx) = \Delta^{-1}\mathbf{1}_{\{|x|>1\}}\Pi(dx)$, it is clear that when $x > 1$, $M(x) = \Delta F^*(-x)$, and $N(x) = \Delta\{1 - F^*(x)\}$. From this one easily checks that

$$A^*(x + \gamma/\Delta) = \frac{1}{\Delta}\int_1^x D(y)\, dy + \int_{-\gamma/\Delta}^1 P(J_1 > y)\, dy$$
$$- \int_{\gamma/\Delta}^1 P(J_1 < -y)\, dy - \int_x^{x+2\gamma/\Delta} P(J_1 < -y)\, dy$$
$$= \frac{A(x) + C}{\Delta} + O(M(x)) \quad \text{as } x \to \infty,$$

where the constant $C$ is given by

$$C = \Delta\left(\int_{-\gamma/\Delta}^1 P(J_1 > y)\, dy - \int_{\gamma/\Delta}^1 P(J_1 < -y)\, dy\right) - \gamma - D(1).$$

A straightforward calculation, treating the cases $|\gamma/\Delta| < 1$, $\gamma \geq \Delta$ and $\gamma \leq -\Delta$ separately, shows that in fact $C = 0$. [For example, in the first case we have $\Delta P(J_1 > y) = N(1)$ and $\Delta P(J_1 < -y) = M(1)$ for $|y| < 1$ so that

$$C = N(1)(1 + \gamma/\Delta) - M(1)(1 - \gamma/\Delta) - \gamma - N(1) + M(1) = 0.]$$

Since $\frac{M(x)}{\sqrt{U(x)M(x)}} = \sqrt{\frac{M(x)}{U(x)}} \to 0$ in all cases, and it is also clear that $\Delta U^*(x) \sim U(x)$ if $U(\infty) = \infty$, and otherwise both $U(\infty)$ and $U^*(\infty)$ are finite, the equivalence of (1.15) and (1.19) follows.

For the other implications, note first that it is known that for any Lévy process $X$ and any fixed $K$, $P(0 \leq X_t \leq K) \to 0$; see, for example, Lemma 2.5 of [4]; thus (1.16) and (1.17) are equivalent. Next, we write

(2.5) $$X_t = S_{N_t}^* + \tilde{X}_t - \gamma t,$$

where of course $S^*$, $\tilde{X}$ and $(N_t, t \geq 0)$ are independent. Since $E\tilde{X}_1 = \gamma$ and $\operatorname{Var} \tilde{X}_1 < \infty$, we know that $P(\tilde{X}_t - \gamma t > 0) \to 1/2$, and it then follows that (1.17) holds if and only if $S^*_{N_t} \xrightarrow{P} +\infty$, and this is easily seen to hold if and only if $S^*_n \xrightarrow{P} +\infty$. (A proof of this statement is given in Lemma 5.2 of [2].) Similarly, (1.18) can be seen to be equivalent to $S^*_n/b_n \xrightarrow{P} +\infty$, and of course these are both equivalent to $P(S^*_n > 0) \to 1$, and hence to (1.15). $\square$

## REFERENCES

[1] BERTOIN, J. (1996). *An Introduction to Lévy Processes*. Cambridge Univ. Press. MR1406564
[2] DONEY, R. A. and MALLER, R. A. (2002). Stability and attraction to Normality for Lévy processes at zero and infinity. *J. Theoret. Probab.* **15** 751–792. MR1922446
[3] DONEY, R. A. and MALLER, R. A. (2002). Moments of passage times of transient Lévy processes. Preprint.
[4] GETOOR, R. K. and SHARPE, M. J. (1994). On the arc-sine law for Lévy processes. *J. Appl. Probab.* **31** 76–89. MR1260572
[5] JANSON, S. (1986). Moments for first-passage and last-exit times, the minimum, and related quantities for random walks with positive drift. *Adv. in Appl. Probab.* **18** 865–879. MR867090
[6] KESTEN, H. (1970). The limit points of a normalized random walk. *Ann. Math. Statist.* **41** 1173–1205. MR266315
[7] KESTEN, H. and MALLER, R. A. (1994). Infinite limits and infinite limit points for random walks and trimmed sums. *Ann. Probab.* **22** 1475–1513. MR1303651
[8] KESTEN, H. and MALLER, R. A. (1996). Two renewal theorems for general random walks tending to infinity. *Probab. Theory Related Fields* **106** 1–38. MR1408415
[9] KESTEN, H. and MALLER, R. A. (1997). Divergence of a random walk through deterministic and random subsequences. *J. Theoret. Probab.* **10** 395–427. MR1455151
[10] KESTEN, H. and MALLER, R. A. (1998). Random walks crossing high level curved boundaries. *J. Theoret. Probab.* **11** 1019–1074. MR1660924
[11] KESTEN, H. and MALLER, R. A. (1998). Random walks crossing power law boundaries. *Studia Sci. Math. Hungar.* **34** 219–252. MR1645198

DEPARTMENT OF MATHEMATICS
UNIVERSITY OF MANCHESTER
MANCHESTER M13 9PL
UNITED KINGDOM
E-MAIL: rad@ma.man.ac.uk